\documentclass[12pt,leqno]{article}
\usepackage{amsfonts}
\usepackage{amssymb}
\usepackage[all]{xy}

\newcommand{\cG}{{\cal G}}

\newcommand{\cM}{{\cal M}}

\newcommand{\cT}{{\cal T}}
\newcommand{\cU}{{\cal U}}

\newcommand{\cX}{{\cal X}}

\newcommand{\Mor}{{\rm Mor}}

\newcommand{\Aut}{{\rm Aut}}

\newcommand{\Spec}{{\rm Spec}}

\newcommand{\ZZ}{{\mathbb Z}}
\newcommand{\CC}{{\mathbb C}}

\newcommand{\RR}{{\mathbb R}}
\newcommand{\NN}{{\mathbb N}}

\newcommand{\PP}{{\mathbb P}}
\newcommand{\QQ}{{\mathbb Q}}

\renewcommand{\wp}{{\mathfrak p}}

\newcommand{\ra}{\rightarrow}
\def\rightepi{{\longrightarrow \kern-0.7em \rightarrow}}

\newcommand{\notteilt}{{\,\not{\kern-0.075em|}\,}}
\def\antiddots{\mathinner{\mkern1mu\raise1pt\vbox{\kern7pt\hbox{.}}\mkern2mu
    \raise4pt\hbox{.}\mkern2mu\raise7pt\hbox{.}\mkern1mu}}
\usepackage{geometry}                
\usepackage{graphicx}
\usepackage{amssymb}
\usepackage{epstopdf}
\begin{document}

\vspace*{15ex}

\begin{center}
{\LARGE\bf Real Belyi Theory} \\
\bigskip
by\\
\bigskip
{\sc  Bernhard K\"ock} and {\sc David Singerman}
\end{center}

\bigskip

\section*{Introduction} Belyi's Theorem \cite{Be} of 1979 had a
profound effect on Galois Theory, Riemann surfaces and complex
algebraic curves. It led Grothendieck  \cite{Gr} to develop his
theory of dessins d'enfants in which there has been a great
interest. The theory is about embedding graphs into compact
Riemann surfaces. For further applications to topics such as such
as moduli spaces and Physics the reader is recommended to consult
\cite{LZ}. The combinatorial and topological foundations of the
theory, as well as the connection with Riemann surfaces, were
developed at the same time by G.~A.~ Jones and D.~Singerman in
\cite{JS1}.

The main purpose of this paper is to develop a Belyi type theory
that applies to Klein surfaces, i.e.\ (possibly non-orientable)
surfaces with boundary which carry a dianalytic structure. More
precisely, let $S$ be a connected compact Klein surface. We prove
that the following statements are equivalent with the possible
exception that statement~(i) may not be implied by the other
statements if $S$ both has empty boundary and is non-orientable:

{\em (i) $S$ can be defined over $ \bar{\QQ} \cap \RR$.

(ii) $S$ admits a Belyi map $\beta: S \ra \Delta$ where $\Delta$
denotes the compactifed closed upper half plane.

(iii) $S$ is isomorphic to the compactification $\overline{\cU/L}$
of the quotient surface $\cU/L$ for some subgroup $L$ of finite
index in the extended modular group $\Gamma^*$ acting on the upper
half plane $\cU$.

(iv) $S$ is isomorphic to the compactification $\overline{\cU/M}$
for some subgroup $M$ of finite index in $\Gamma^*(2)$, the
extended principal congruence subgroup of level 2.

(v) $S$ is isomorphic to the quotient surface $\cX/K$ for some
subgroup $K$ of finite index in an extended triangular group
acting on one of the three simply connected Riemann surfaces $\cX=
{\hat{\CC}}$, $\CC$ or $\cU$.

(vi) $S$ carries a map $\cM$ such that $S$ is isomorphic to the
quotient surface $\cX/K$ where $K$ is the map subgroup for $\cM$.}

We refer the reader to the main body of this paper for a detailed
definition of the various notions which have been used in the
formulation of these statements. The equivalence between (i) and
(ii), i.e.\ the real analogue of Belyi's theorem, is discussed in
Section~2. In Section~3 we represent Klein surfaces as quotient
surfaces. And finally, in Section~4, we develop a theory of maps
on Klein surfaces; the combinatorial and topological foundations
have been laid out by R.~P.~Bryant and D.~Singerman in \cite{BS}.

{\em Acknowledgements.} We would like to thank Gareth Jones,
Javier Cirre and Pablo Mart\'in for many helpful ideas.

\bigskip

\bigskip

\section*{\S 1 Belyi's Theorem for curves over non algebraically
closed base fields}

In this section we extend Belyi's theorem from complex curves to
curves over an arbitrary subfield of $\CC$. We follow the approach
used in \cite{Ko} and based on earlier work by Wolfart, see
\cite{Wo}. Related results may be found in the recent paper
\cite{Go} by Gonz\'alez-Diez.

We fix a subfield $C$ of $\CC$. In this section, by a {\em curve
over $C$} we mean a smooth projective geometrically connected
variety of dimension~1 over $C$. We recall that mapping a curve
$X$ over $C$ to its function field $K(X)$ yields an
antiequivalence between the category of curves over $C$ and
non-constant morphisms on the one hand and the category of
finitely generated field extensions $K$ of $C$ of transcendence
degree~1 such that $C$ is algebraically closed in $K$ on the other
hand. Furthermore we will frequently use the well-known fact that
elements of $K(X)$ correspond to morphisms $X \ra \PP^1_C$ of
varieties over $C$. As usual we say that a curve $X$ over $C$ {\em
can be defined over a subfield $D$ of $C$} if there is a curve $Y$
over $D$ such that $Y_C := Y \times_D C$ is isomorphic to $X$ over
$C$. This is obviously equivalent to the condition that there is a
subfield $K$ of the function field $K(X)$ containing $D$ which is
finitely generated and of transcendence degree~1 over $D$ such
that $D$ is algebraically closed in $K$ and such that the
canonical homomorphism
\[{\rm Quot}(C\otimes_D K) \ra K(X)\]
from the field of fractions ${\rm Quot}(C \otimes_D K)$ of $C
\otimes_D K$ to $K(X)$ is bijective. Note that $C\otimes_D K$ is
an integral domain by the following lemma.

{\bf (1.1) Lemma.} {\em Let $D$ be any field of characteristic~0
and let $K$ be any extension field of $D$. Then $D$ is
algebraically closed in $K$ if and only if, for every extension
field $C$ of $D$, the tensor product $C\otimes_D K$ is an integral
domain. In this case, $C$ is algebraically closed in ${\rm Quot}(C
\otimes_D K)$.}

{\em Proof.} See Proposition~(4.3.2) and Proposition~(4.3.5) in
\cite{EGAIV}. \hspace*{\fill} $\Box$

We recall (see \S 2 in \cite{Ko}) that the {\em moduli field
$M(X,t)$} of a finite morphism $t: X \ra \PP^1_\CC$ from a curve
$X$ over $\CC$ to the complex projective line $\PP^1_\CC$ is
defined to be the subfield of $\CC$ fixed by the subgroup $U(X,t)$
of $\Aut(\CC)$ consisting of all automorphisms ${\sigma}$ of $\CC$
such that there is an isomorphism $f_{\sigma}: X^{\sigma} \ra X$
of varieties over $\CC$ such that the following diagram commutes:
\[ \xymatrix{X^{\sigma} \ar[rr]^{f_{\sigma}} \ar[d]^{t^{\sigma}} && X \ar[d]^t \\
(\PP^1_\CC)^{\sigma} \ar[rr]^{\Spec({\sigma})} && \PP^1_\CC}\]
Here, the notation $X^{\sigma}$ means the scheme $X$ viewed as a
variety over $\CC$ via the new structure morphism $X \,\,
\stackrel{p}{\ra} \,\, \Spec(\CC) \,\,
\stackrel{\Spec({\sigma})}{\longrightarrow} \,\, \Spec(\CC)$ where
$p$ denotes the structure morphism of the given variety $X$ over
$\CC$ (see \S 1 in \cite{Ko} for further explanations of this
notation).

The following proposition is a special case of Corollary~(3.2) in
\cite{Ko}.

{\bf (1.2) Proposition.} {\em If the critical values of $t$ are
$\QQ$-rational then the moduli field $M(X,t)$ is a number field.}
\hspace*{\fill} $\Box$

Let now $t: X \ra \PP^1_C$ be a finite morphism from a curve $X$
over the given subfield $C$ of $\CC$ to the projective line
$\PP^1_C$. Then the subgroup $U(X_\CC,t_\CC)$ of $\Aut(\CC)$
obviously contains $\Aut(\CC/C)$. Thus the moduli field
$M(X_\CC,t_\CC)$ of $t_\CC: X_\CC \ra \PP^1_\CC$ is contained in
$\CC^{\Aut(\CC/C)} = C$ (by Lemma~(1.4) in \cite{Ko}).

{\bf (1.3) Proposition.} {\em If there is a $C$-rational
unramified point $P$ on $X$ such that $Q:= t(P)$ lies in $\QQ$
then the curve $X/C$ and the morphism $t$ can be defined over a
finite extension of the moduli field $M(X_\CC, t_\CC)$ (inside
$C$).}

Note that Proposition~(1.3) is stronger than Theorem~(2.2) in
\cite{Ko} which only implies that $X_\CC/\CC$ (or
$X_{\bar{C}}/\bar{C}$) is defined over a finite extension of
$M(X_\CC, t_\CC)$. For the convenience of the reader we give the
details of the proof of Proposition~(1.3) which is quite similar
to the proof of Theorem~(2.2) in \cite{Ko}.

{\em Proof.} By the Riemann-Roch Theorem (see \cite{Mo}) applied
to the divisor $D:= ({\rm genus}(X) + 1) [P]$ there is a
non-constant meromorphic function $z$ on $X$ (defined over $C$)
such that $P$ is the only pole of $z$. Then the function field
$K(X)$ of $X$ is generated over $C$ by $t$ and $z$. This follows
from the easy observation that the field extension $K(X)/C(t,z)$
is a subextension of $K(X)/C(t)$ and of $K(X)/C(z)$ and that,
hence, the corresponding morphism of curves over $C$ is both
unramified and totally ramified at $P$. We assume furthermore that
we have chosen $z$ in such a way that the pole order $m:= -{\rm
ord}_P(z) \in \NN$ is minimal. Then the constant function $1$ and
the function $z$ form a basis of the $C$-vector space
\[V:= \{ x \in K(X): x\textrm{ has a pole at most at $P$ and }
{\rm ord}_P(x) \ge -m\}.\] To see this, let $x \in V$; if $x$ is
not a constant we have ${\rm ord}_P(x) = -m \,\, \stackrel{{\rm
def}}{=}\,\, {\rm ord}_P(z)$ since $m$ was minimal; hence there is
a constant $\alpha \in C$ such that $-{\rm ord}_P(x-\alpha z) <
m$; thus, again since $m$ was minimal, the function $x- \alpha z$
is constant, as desired.\\
Since $t$ is unramified at $P$, the meromorphic function $t-Q$ on
$X$ is a local parameter on $X$ at $P$. Obviously, there is a
unique function $z' \in V$ such that the leading coefficient
(i.e., the coefficient of $(t-Q)^{-m}$) and the constant
coefficient (i.e., the coefficient of $(t-Q)^0$) in the Laurent
expansion of $z'$ with respect to the local parameter $t-Q$ are
$1$ and $0$, respectively. We may and we will assume that
$z=z'$.\\
We now claim that the minimal polynomial of $z$ over $C(t)$ has
coefficients in $k(t)$ where $k$ is a finite extension of
$M(X_\CC/t_\CC)$ inside $C$. Then the field extension $K(X)/C(t)$
is defined over $k$. This means that $X/C$ and the morphism $t$
are defined over $k$, as was to be shown.\\
To prove the above claim we introduce the notation
$U(X_\CC,t_\CC,P)$ for the subgroup of $U(X_\CC, t_\CC)$
consisting of all automorphisms ${\sigma}$ of $\CC$ such that
there is an isomorphism $f_{\sigma}: X_\CC^{\sigma} \ra X_\CC$ of
varieties over $\CC$ such that the diagram
\[\xymatrix{X_\CC^{\sigma} \ar[rr]^{f_{\sigma}} \ar[d]^{t_\CC^{\sigma}}
&& X_\CC \ar[d]^{t_\CC} \\ (\PP^1_\CC)^{{\sigma}}
\ar[rr]^{\Spec({\sigma})} && \PP^1_\CC }\] commutes and such that
$f_{\sigma}(P^{\sigma}) = P$; here, $P^{\sigma}$ denotes the point
in $X_\CC^{\sigma}$ corresponding to $P$. Note that $f_{\sigma}$
is unique since $\Aut(t_\CC)$ acts freely on the fibre
$t_\CC^{-1}(Q)$. Thus mapping ${\sigma}$ to the automorphism of
the function field $K(X_\CC)$ induced by $f_{\sigma}$ yields an
action of $U(X_\CC,t_\CC,P)$ on $K(X_\CC)$ by semilinear field
automorphisms which fix $t_\CC \in K(X_\CC)$. The meromorphic
function $z_\CC \in K(X_\CC)$ is invariant under the action of
$U(X_\CC, t_\CC,P)$ defined above since the image of $z_\CC$ under
${\sigma} \in U(X_\CC, t_\CC,P)$ has the same three defining
properties as $z_\CC$, as one easily checks. Hence the minimal
polynomial of $z_\CC$ over $\CC(t_\CC)$ is invariant under the
action of $U(X_\CC, t_\CC,P)$, i.e.\ it has coefficients in
$k(t_\CC)$ where $k$ is the subfield of $\CC$ fixed by $U(X_\CC,
t_\CC,P)$. We obviously have $\Aut(\CC/C) \subseteq U(X_\CC,
t_\CC, P)$ since $P$ is $C$-rational. Thus $k$ is contained in
$\CC^{\Aut(\CC/C)} = C$ (by Lemma~(1.4) in \cite{Ko}). Finally,
being the stabilizer of $[P]$ under the (well-defined!) action
$({\sigma}, [P]) \mapsto [f_{\sigma}(P)]$ of $U(X_\CC, t_\CC)$ on
$t^{-1}_\CC(Q)/\Aut(t_\CC)$, the subgroup $U(X_\CC, t_\CC, P)$ has
finite index in $U(X_\CC, t_\CC)$. Hence, by Lemma~(1.6) in
\cite{Ko}, $k$ is a finite extension of $M(X_\CC, t_\CC)$, as was
to be shown. \hspace*{\fill} $\Box$

{\bf (1.4) Theorem (Belyi's theorem for curves over a non
algebraically closed base field)}. {\em Lex $X$ be a curve over
$C$.
\\
(a) If $X$ can be defined over $\bar{\QQ} \cap C$, then there
exists a finite morphism $t:X \ra \PP^1_C$ (defined over $C$!)
such that the critical values of $t_\CC: X_\CC \ra \PP^1_\CC$ lie
in $\{0,1,\infty\}$. (The morphism $t$ may actually be chosen to
be defined even over $\bar{\QQ} \cap C$.)\\
(b) The converse of statement~(a) is true if there exists an
unramified $C$-rational point $P$ on $X$ such that $t(P)$ lies in
$\QQ$.}

{\em Proof.} \\
(a) Let $X_{\bar{\QQ}\cap C}$ be a model of $X$ over
$\bar{\QQ}\cap C$ and let $t': X_{\bar{\QQ}\cap C} \ra
\PP^1_{\bar{\QQ} \cap C}$ be any finite morphism. By Lemma~(3.4)
in \cite{Ko}, the critical values of $t'_\CC: X_\CC \ra \PP^1_\CC$
lie in $\bar{\QQ} \cup \{\infty\}$. By Lemma~(3.5) and (the proof
of) Lemma~(3.6) in \cite{Ko} there is a finite morphism
$r:\PP^1_\CC \ra \PP^1_\CC$ defined over $\QQ$ such that the
critical values of $r$ and the images of the critical values of
$t'_\CC$ under $r$ lie in $\{0,1,\infty\}$. We write $r$ also for
the corresponding morphism $r: \PP^1_{\bar{\QQ}\cap C} \ra
\PP^1_{\bar{\QQ}\cap C}$. Let $t$ denote the composition $r \circ
t'$. Then $t_C: X \ra \PP^1_C$ is a finite morphism defined over
$\bar{\QQ} \cap C$ such that the
critical values of $t_\CC$ lie in $\{0,1, \infty\}$. \\
(b) By Proposition~(1.2), the moduli field $M(X_\CC, t_\CC)$ of
$t_\CC$ is a number field (contained in $C$). By
Proposition~(1.3), the curve $X$ and the morphism $t:X \ra
\PP^1_\CC$ can be defined over a finite extension of $M(X_\CC,
t_\CC)$ (inside $C$). Thus $X$ and $t$ can be defined over
$\bar{\QQ} \cap C$, as was to be shown. \hspace*{\fill} $\Box$

\bigskip

\section*{\S 2 Belyi's Theorem for Klein surfaces}

In this section we develop a `real' Belyi theory; its highlight is
a Belyi theorem for Klein surfaces, see Theorem~(2.6).

We recall that the canonical functor from the category of smooth
projective curves over $\CC$ to the category of compact Riemann
surfaces is an equivalence of categories. In particular, every
compact connected Riemann surface $X$ is isomorphic to the Riemann
surface associated with a smooth projective connected curve over
$\CC$ and the field $M_{{\rm RS}}(X)$ of meromorphic functions on
$X$ can be identified with the function field of this curve over
$\CC$.

As usual we say that a compact Riemann surface $X$ {\em can be
defined over a subfield $C$ of $\CC$} if the corresponding smooth
projective connected curve over $\CC$ can be defined over $C$. A
reader not familiar with the language of schemes may prefer to
work with the following characterization. By \S 1 and the remarks
after Definition~(4.5.2) in \cite{EGAIV}, the Riemann surface $X$
can be defined over $C$ if and only if there exists a subfield $K$
of $M_{{\rm RS}}(X)$ containing $C$ which is finitely generated
and of transcendence degree~1 over $C$ such that the canonical map
$\CC \otimes_C K \ra M_{{\rm RS}}(X)$ is injective and $M_{{\rm
RS}}(X)$ is the field of fractions of $\CC \otimes_C K$ via this
map.

Let now $S$ be a compact connected Klein surface. We refer the
reader to \cite{AG} and \cite{Ga} for basic facts about Klein
surfaces. We denote the field of meromorphic functions of $S$ by
$M_{{\rm KS}}(S)$. By Theorem~1.4.6 in \cite{AG}, this is the set
of morphisms (of Klein surfaces) from $S$ to the compactified
closed upper half plane
\[\Delta:= \{z \in \CC: {\rm Im}(z) \ge 0\} \cup
\{\infty\}.\] $M_{\rm KS}(S)$ is a finitely generated extension of
$\RR$ of transcendence degree 1.

{\bf (2.1) Definition.} We say that {\em $S$ can be defined over
the field $\bar{\QQ} \cap \RR$} of real algebraic numbers if there
exists a subfield $K$ of $M_{{\rm KS}}(S)$ containing $\bar{\QQ}
\cap \RR$ which is finitely generated and of transcendence
degree~1 over $\bar{\QQ} \cap \RR$ such that the canonical map
$\RR \otimes_{\bar{\QQ}\cap \RR} K \ra M_{{\rm KS}}(S)$ is
injective and $M_{{\rm KS}}(S)$ is the field of fractions of $\RR
\otimes_{\bar{\QQ}\cap \RR} K$ via this map.

Note that $S$ can obviously be defined over $\bar{\QQ} \cap \RR$
if and only if the projective real algebraic curve corresponding
to $M_{\rm KS}(S)$ can be defined over $\bar{\QQ} \cap \RR$. When
using this characterization, a real projective algebraic curve has
basically to be understood as a scheme, i.e.\ as a space {\em
together} with the defining equations. In contrast to complex
curves, it is not sufficient to think of a real curve just as the
set of real points, nor as the set of its complex points. In order
to make this section (and this paper) accessible to as wide an
audience as possible we have decided to work with the language of
function fields.

{\bf (2.2) Lemma.} {\em Let $X$ be a compact connected Riemann
surface. Then $X$ can be defined over $\bar{\QQ}$ (as Riemann
surface) if and only if $X$, viewed as Klein surface, can be
defined over $\bar{\QQ} \cap \RR$.}

{\rm Proof.} We first prove the only-if-part. Let $K$ be a
subfield of $M_{{\rm RS}}(X)$ containing $\bar{\QQ}$ which is
finitely generated and of transcendence degree~1 over $\bar{\QQ}$
such that the canonical map $\CC \otimes_{\bar{\QQ}} K \ra M_{{\rm
RS}}(X)$ is injective and $M_{{\rm RS}}(X)$ is the field of
fractions of $\CC \otimes_{\bar{\QQ}} K$. The canonical
$\RR$-algebra homomorphism $\RR \otimes_{\bar{\QQ}\cap\RR}
\bar{\QQ} \ra \CC$ is bijective since it takes the $\RR$-basis
$1\otimes 1, 1 \otimes i$ of $\RR \otimes_{\bar{\QQ}\cap\RR}
\bar{\QQ}$ to the $\RR$-basis $1,i$ of $\CC$. Hence we obtain
\[\RR \otimes_{\bar{\QQ}\cap\RR} K \cong
\RR\otimes_{\bar{\QQ}\cap\RR} \bar{\QQ} \otimes_{\bar{\QQ}} K
\cong \CC \otimes_{\bar{\QQ}} K\] and
\[{\rm Quot}(\RR \otimes_{\bar{\QQ}\cap\RR} K) \cong {\rm
Quot}(\CC \otimes_{\bar{\QQ}} K) \cong M_{{\rm RS}}(X) \cong
M_{{\rm KS}}(X), \] as desired.\\
We now assume that $X$, viewed as Klein surface, can be defined
over $\bar{\QQ} \cap \RR$. Let $K$ be a subfield of $M_{{\rm
KS}}(X)$ containing $\bar{\QQ}\cap\RR$ which is finitely generated
and of transcendence degree~1 over $\bar{\QQ}\cap\RR$ such that
the canonical map $\RR\otimes_{\bar{\QQ}\cap\RR} K \ra M_{{\rm
KS}}(X)$ is injective and such that $M_{{\rm KS}}(X)$ is the field
of fractions of $\RR \otimes_{\bar{\QQ}\cap\RR} K$. We have
$i=\sqrt{-1} \in K$ because otherwise the ring
\[\RR \otimes_{\bar{\QQ}\cap\RR} K \otimes_\QQ \QQ[i] \cong \RR
\otimes_{\bar{\QQ}\cap\RR} K[i] \cong \RR
\otimes_{\bar{\QQ}\cap\RR} \bar{\QQ} \otimes_{\bar{\QQ}} K[i]
\cong \CC \otimes_{\bar{\QQ}} K[i]\] would be an integral domain
by Lemma~(1.1); but the canonical map
\[{\rm Quot}(\RR \otimes_{\bar{\QQ}\cap\RR} K) \otimes_\QQ \QQ[i]
\ra {\rm Quot}(\RR \otimes_{\bar{\QQ}\cap\RR} K \otimes_\QQ
\QQ[i])\] is injective (We may use the basis $1,i$ of $\QQ[i]$
over $\QQ$ to check this.) and the ring
\[{\rm Quot}(\RR \otimes_{\bar{\QQ}\cap\RR} K) \otimes_\QQ \QQ[i]
\cong M_{{\rm KS}}(X) \otimes_\QQ \QQ[i] \cong M_{{\rm RS}}(X)
\otimes_{\QQ} \QQ[i] \cong M_{{\rm RS}}(X) \times M_{{\rm
RS}}(X)\] is not an integral domain. Hence $K$ contains $\bar{\QQ}
= (\bar{\QQ}\cap\RR)[i]$ and we have
\[{\rm Quot}(\CC \otimes_{\bar{\QQ}} K) \cong {\rm Quot}(\RR
\otimes_{\bar{\QQ}\cap\RR} K) \cong M_{{\rm KS}}(X) \cong M_{{\rm
RS}}(X),\] as desired. \hspace*{\fill} $\Box$

We recall that the canonical functor from the category of Riemann
surfaces to the category of Klein surfaces admits a right adjoint
functor. In other words, for every Klein surface $S$ there exists
a Riemann surface $S^{{\rm c}}$ together with a morphism $f_S:
S^{{\rm c}} \ra S$ of Klein surfaces  such that for every Riemann
surface $X$ the map
\[\Mor_{{\rm RS}}(X,S^{{\rm c}}) \ra \Mor_{{\rm KS}}(X,S), \quad
\rho \mapsto f_S \circ \rho,\] is bijective (see Proposition~1.6.2
in \cite{AG}). In fact the map $f_S$ is a double cover; we
therefore call $S^{\rm c}$ the {\em complex double of $S$}. The
complex double $S^{\rm c}$ comes with an antiholomorphic
involution $\tau$ such that $S^{\rm c}/\langle \tau \rangle \cong
S$. In fact this gives an equivalence between the category of
Klein surfaces and the category of pairs consisting of a Riemann
surface and an antiholomorphic involution. For instance,
$\Delta^{{\rm c}}$ is the Riemann sphere ${\hat{\CC}}$ and
$f_\Delta: \Delta^{{\rm c}} \ra \Delta$ is the {\em folding map}
${\hat{\CC}} \ra \Delta, a+bi \mapsto a +|b|i, \infty \mapsto
\infty$; or, if $S$ is the Klein surface associated with a Riemann
surface $X$, then $S^{{\rm c}}$ is the disjoint union of the
Riemann surfaces $X$ and $\bar{X}$ (where $\bar{X}$ denotes the
topological space $X$ equipped with the antiholomorphic structure)
and $f_S: X \cup \bar{X} \ra S$ is the identity map on both $X$
and $\bar{X}$.

{\bf (2.3) Definition.} A morphism $\beta: S \ra \Delta$ from a
compact Klein surface $S$ to the compactified upper half plane
$\Delta$ is called a {\em Belyi map on $S$} if the complex double
$\beta^{{\rm c}}: S^{{\rm c}} \ra \Delta^{{\rm c}} \cong
{\hat{\CC}}$ of $\beta$ is a Belyi map on the Riemann surface
$S^{{\rm c}}$, i.e. if the restriction of $\beta^{{\rm c}}$ to
every connected component of $S^{{\rm c}}$ has at most 3 critical
values. We shall also call $\beta$  {\em a real Belyi function.}

{\bf (2.4) Lemma.} {\em Let $X$ be a compact connected Riemann
surface. Then $X$ (as Riemann surface) admits a Belyi map from $X$
to ${\hat{\CC}}$ if and only if $X$, viewed as Klein surface,
admits a Belyi map from $X$ to $\Delta$.}

{\em Proof.} We first prove the if-part. Let $\beta: X \ra \Delta$
be a Belyi map on the Klein surface $X$. Then, by definition, the
induced morphism $\beta^{{\rm c}}: X^{{\rm c}} \ra \Delta^{{\rm
c}} \cong {\hat{\CC}}$ is a Belyi map on the Riemann surface
$X^{{\rm c}}$. In particular, the restriction of $\beta^{{\rm c}}$
to the connected component $X$ of $X^{{\rm c}} = X \cup \bar{X}$
is a Belyi map on the given Riemann surface $X$, as desired.\\
We now prove the only-if-part. Let $\beta: X \ra {\hat{\CC}}$ be a
Belyi map on the Riemann surface~$X$. Then the composition $X \,\,
\stackrel{\beta}{\ra} \,\, {\hat{\CC}} \,\, \stackrel{\rho}{\ra}
\,\, \Delta$ of $\beta$ with the folding map $\rho: {\hat{\CC}}
\ra \Delta$ is a Belyi map on the Klein surface $X$ since the
induced morphism $(\rho \circ \beta)^{{\rm c}}: X \cup \bar{X} =
X^{{\rm c}} \ra \Delta^{{\rm c}} \cong {\hat{\CC}}$ is equal to
$\beta$ on $X$ and equal to the complex conjugate $\bar{\beta}$ of
$\beta$ on $\bar{X}$ and, hence, is a Belyi map on the Riemann
surface $X^{{\rm c}}$. \hspace*{\fill} $\Box$

{\bf (2.5) Lemma.} {\em If a compact Klein surface $S$ admits a
Belyi map $\beta: S \ra \Delta$ then it also admits a Belyi map
$\gamma: S\ra \Delta$ such that the critical values of
$\gamma^{{\rm c}}: S^{{\rm c}} \ra {\hat{\CC}}$ lie in
$\{0,1,\infty\}$. }

{\em Proof.} We may assume that $S$ is connected. \\
We first consider the case that the critical values of
$\beta^{{\rm c}}: S^{{\rm c}} \ra {\hat{\CC}}$ lie in $\RR \cup
\{\infty\}$. Then we can find a M\"obius transformation
$r:{\hat{\CC}} \ra {\hat{\CC}}$ with real coefficients which maps
the critical values of $\beta^{{\rm c}}$ into $\{0,1,\infty\}$. In
particular, the unique $\gamma \in M_{{\rm KS}}(S)$ such that
$\gamma^{{\rm c}} = r \circ \beta^{{\rm c}}$ is a Belyi map on $S$
of the desired type. \\
If a critical value of $\beta^{{\rm c}}$, say $P$, lies in $\CC
\backslash \RR$ then also its complex conjugate $\bar{P}$ is a
critical value of $\beta^{{\rm c}}$ since $\beta^{{\rm c}}$ is
defined over $\RR$. By the same argument, the third one (if there
is one), say $Q$, must lie in $\RR \cup \{\infty\}$. As above, by
composing with a M\"obius transformation with real coefficients,
we may assume that $Q = \infty$. Let $a:= -2{\rm Re}(P) \in \RR$
and $b:= |P|^2 \in \RR$. Then the real quadratic polynomial
$r:=X^2+aX+b$ takes the three critical values $P,\bar{P},Q=\infty$
of $\beta^{{\rm c}}$ to the two values $0,\infty$ and it has the
two critical values $b-\frac{a^2}{4}$ and $\infty$. Thus the
unique $\tilde{\beta} \in M_{{\rm KS}}(S)$ such that
$\tilde{\beta}^{{\rm c}} = r \circ \beta^{{\rm c}}$ is a Belyi map
on $S$ of the type considered in the first case. Therefore $S$
admits a Belyi map of the desired type also in the general case.
\hspace*{\fill} $\Box$

{\bf (2.6) Theorem (Belyi's Theorem for Klein surfaces).} {\em Let
$S$ be a compact connected Klein surface. If $S$ can be defined
over $\bar{\QQ} \cap \RR$ then $S$ admits a Belyi map $\beta: S
\ra \Delta$. If the boundary $\partial S$ of $S$ is non-empty or
if $S$ is the Klein surface associated with a compact connected
Riemann surface then the converse is true as well.}

{\em Proof.} If $S$ is the Klein surface associated with a Riemann
surface then Theorem~(2.6) follows from Lemma~(2.2), Lemma~(2.4)
and Belyi's Theorem for Riemann surfaces (see Theorem~(3.3) in
\cite{Ko}).\\
Hence we may assume that $S$ is not the Klein surface associated
with a Riemann surface. Then its complex double $S^{{\rm c}}$ is
connected by Lemma~1.6.3 in \cite{AG} and Remark~(1.5)(ii) in
\cite{Ga}; in particular $M_{{\rm RS}}(S^{{\rm c}})$ is a field
and we have $M_{{\rm RS}}(S^{{\rm c}}) = M_{{\rm KS}}(S)
\otimes_\QQ \QQ[i]$ by Corollary~1.6.5 in \cite{AG}. In other
words, the smooth projective curve $X$ over $\RR$ corresponding to
the function field $M_{{\rm KS}}(S)$ over $\RR$ is geometrically
connected, and we may apply the results of \S 1.\\
If $S$ can be defined over $\bar{\QQ} \cap \RR$ then obviously the
curve $X$ can be defined over $\bar{\QQ} \cap \RR$. Hence, by
Theorem~(1.4)(a), there exists a finite morphism $t:X \ra
\PP^1_\RR$ (defined over $\RR$) such that the critical values of
$t_\CC: X_\CC \ra \PP^1_\CC$ lie in $\{0,1,\infty\}$. Thus the
corresponding morphism $\beta: S \ra \Delta$ is a Belyi map on the
Klein surface $S$, as desired.\\
We now assume that $S$ admits a Belyi map $\beta: S \ra \Delta$
and that $\partial S \not= \emptyset$. As $\partial S$ is
homoemorphic to the real manifold $X(\RR)$ of $\RR$-valued points
on $X$ (see Theorem~(4.17) in \cite{Ga}) there exists an
$\RR$-valued point on $X$, say $P_0$. Since the manifold $X(\RR)$
is of dimension~1 there exists a neighborhood $U$ of $P_0$ in
$X(\RR)$ homeomorphic to an interval in $\RR$. The morphism $t: X
\ra \PP^1_\RR$ corresponding to the given Belyi map $\beta: S \ra
\Delta$ maps $U$ to an interval in $\RR \cup \{\infty\}$. In
particular there are infinitely many points $P$ in $U \subseteq
X(\RR)$ such that $t(P)$ lies in $\QQ$. Since the set of ramified
points is finite, there hence exists an unramified $\RR$-rational
point $P$ on $X$ such that $t(P)$ lies in $\QQ$. Thus, by
Theorem~(1.4)(b), the curve $X$ and hence the Klein surface $S$
can be defined over $\bar{\QQ}\cap\RR$, as was to be shown.
\hspace*{\fill} $\Box$

{\bf (2.7) Question.} Let $S$ be a compact connected Klein surface
which admits a Belyi map. By Theorem~(2.6), the Klein surface $S$
can be defined over $\bar{\QQ}\cap\RR$ if $\partial S \not=
\emptyset$ or if ($\partial S = \emptyset$ and $S$ is orientable).
Does the same conclusion hold also in the remaining case that
$\partial S = \emptyset$ and $S$ is not orientable?

{\bf (2.8) Example.} Let $S$ be a compact connected Klein surface
which is non-orientable or has non-empty boundary. Then the
complex double $S^{\rm c}$ is connected again. Let $g$ denote the
genus of the compact Riemann surface $S^{\rm c}$. In this example
we want to show that Question~(2.7) has an affirmative answer if
$g=0$ or $g=1$. \\
(a) Let $g = 0$, i.e.\ the Riemann surface $S^{\rm c}$ is
isomorphic to the Riemann sphere~${\hat{\CC}}$. Then, according to
Theorem 1.9.4 in \cite{AG}, the Klein surface $S$ is isomorphic to
either the compactified upper halfplane $\Delta$ or to the real
projective plane $\PP^2(\RR)$ which, as Klein surface, is defined
as ${\hat{\CC}}/\langle {\sigma} \rangle$ where the
antiholomorphic involution $\sigma$ sends $z \in {\hat{\CC}}$ to
$-\frac{1}{\bar{z}}$. In the first case, the function field of $S$
is given by the field $\RR(z)$ of rational functions in one
variable. In the second case, the function field of $S$ is
isomorphic to the field of fractions of the integral domain
$\RR[X,Y]/(X^2+Y^2+1)$ embedded into the function field $\CC(z)$
of the Riemann sphere via $X \mapsto \frac{1}{2}\left(z -
\frac{1}{z}\right)$ and $Y \mapsto \frac{1}{2i}\left(z +
\frac{1}{z}\right)$. Both of these Klein surfaces are obviously
defined over $\QQ$, the smallest field of characteristic $0$. In
particular, Question~(2.7) has an affirmative answer if $S$ is the
real projective plane. \\
(b) Let $g=1$, i.e.\ the Riemann surface $S^{\rm c}$ is an
elliptic curve. From \S 9 in \cite{AG} we know $S$ that must be
homoemorphic to either an annulus, a M\"obius strip or to a Klein
bottle. By Example~1.7 in \cite{Ko}, the moduli field of $S^{\rm
c}$ is $\QQ(j)$ where $j \in \RR$ denotes the $j$-invariant of
$S^{\rm c}$. According to Chapters~14.34, 14.41 and 14.42 in
\cite{Al}, the function field of $S$ is isomorphic to the field of
fractions of an integral domain of the form $\RR[X,Y]/\left(Y^2
\pm (1 \pm X^2)(1 \pm \lambda X^2)\right)$ where $\lambda \in \RR$
is a/the Legendre modulus of the elliptic curve $S^{\rm c}$. Here
$X$ corresponds to (a constant multiple of) the Jacobi elliptic
function ${\rm sn}$ (`sinam function') and $Y$ to its derivative.
As $\lambda$ satisfies the equation
\[ j = k \frac{(\lambda^2 - \lambda + 1)^3}{\lambda^2 (\lambda -
1)^2} \qquad \textrm{ with some } k \in \QQ\] (see Example~1 on
p.~566 in \cite{JS2}) we see that $S$ can be defined over a (real)
extension field of $\QQ(j)$ of degree at most $6$. Now
Proposition~(1.2) implies that Question~(2.7) has an affirmative
answer also if $S$ is homeomorphic to the Klein bottle.

{\bf (2.9) Example.} Let $S$ and $j$ be as in Example~(2.8)(b)
above. We now want to show that $S$ can be defined not only over a
finite extension of $\QQ(j)$ but over $\QQ(j)$ itself if $S$ is
homeomorphic to an annulus or to a M\"obius strip. Whether the
same holds true if $S$ is homeomorphic to a Klein bottle remains
open.\\
We consider the two real elliptic curves $X_+(j)$ and $X_-(j)$
given by the following equations defined over $\QQ(j)$:
\begin{eqnarray*}
X_+(j): & y^2 = 4x^3 - \frac{27j}{j-1} x - \frac{27j}{j-1} & \textrm{if } j\not= 0,1\\
&         y^2 = 4x^3 - 1 & \textrm{if } j= 0\\
&         y^2 = 4x^3 - x & \textrm{if } j= 1\\
X_-(j): & y^2 = 4x^3 - \frac{27j}{j-1} x + \frac{27j}{j-1} &
\textrm{if } j \not= 0,1\\
&         y^2 = 4x^3 + 1 & \textrm{if } j =0\\
&         y^2 = 4x^3 + x & \textrm{if } j =1
\end{eqnarray*}
It is easy to see that the $j$-invariant of both $X_+(j)$ and
$X_-(j)$ is $j$. However $X_+(j)$ is not $\RR$-isomorphic to
$X_-(j)$. To see this we recall that two real elliptic curves $X$
and $\tilde{X}$ given by the equations $y^2 = 4x^3 - g_2x - g_3$
and $y^2 = 4 x^3 - \tilde{g}_2 x - \tilde{g}_3$ are
$\RR$-isomorphic if and only if there exists a $u \in \RR$ such
that the second equation is obtained from the first equation by
the substitutions $x \mapsto u^2 x$ and $y \mapsto u^3 y$. In
particular, if $X$ and $\tilde{X}$ are isomorphic then the sign of
$g_2$ is equal to the sign of $\tilde{g}_2$ and the sign of $g_3$
is equal to the sign of $\tilde{g}_3$. As this is not true for
$X_+(j)$ and $X_-(j)$, they are not $\RR$-isomorphic. Obviously,
both $X_+(j)$ and $X_-(j)$ have $\RR$-valued points; hence the
corresponding Klein surfaces $S_+(j)$ and $S_-(j)$ have non-empty
boundary and are therefore not homeomorphic to a Klein bottle.
From the classification of Klein surfaces of genus~1 (see \S 9 of
Chapter~1 in \cite{AG}) we know that there are exactly two
isomorphism classes of connected compact Klein surfaces which are
not homeomorphic to a Klein bottle and whose complex double is
connected and has genus~1 and $j$-invariant $j$. Hence $S$ is
isomorphic to $S_+(j)$ or $S_-(j)$ and therefore defined over
$\QQ(j)$.

\bigskip

\bigskip

\bigskip

\section*{\S 3 Belyi Functions on Klein Surfaces and the Extended
Modular Group}

The object of this section is to extend various characterizations
of Belyi surfaces (see Section~4 in \cite{JS2}) from Riemann
surfaces to Klein surfaces.

Let $\Gamma = {\rm PSL}_2(\ZZ)$ denote modular group which we regard
as usual as M\"obius transformations on the complex plane $\CC$ or
upper half plane $\mathcal{U}:= \{z \in \CC : \rm{Im}(z) > 0\}$. Let
$\Gamma^*$ be the group generated by $\Gamma$ and the
antiholomorphic reflection $R: z \mapsto -\bar{z}$. We call
$\Gamma^*$ the {\em extended modular group}. It is easy to see that
$\Gamma$ is of index $2$ in $\Gamma^*$.

{\bf (3.1) Lemma.} {\em The $J$-function $J: \mathcal{U}
\rightarrow \CC$ induces an isomorphism of Klein surfaces
\[\mathcal{U}/{\Gamma^*} \,\, \stackrel{\sim}{\rightarrow} \,\,
\Delta \backslash \{\infty\} .\]}

{\em Proof.} It is well-know that $J$ induces an isomorphism
$\mathcal{U}/\Gamma \rightarrow \CC$. It therefore suffices to show
that the following diagram commutes:
\[\xymatrix{\mathcal{U} \ar[r]^R \ar[d]^J & \mathcal{U} \ar[d]^J \\
\CC \ar[r]^{\rm c} & \CC}\] where $c$ denotes complex conjugation.
This easily follows from the fact that for every $\tau \in \cU$
the lattice $\ZZ+\ZZ(-\bar{\tau})$ is equal to the complex
conjugate of the lattice $\ZZ + \ZZ\tau$. \hfill $\Box$

{\bf (3.2) Theorem.} {\em A compact Klein surface $S$ admits a
Belyi map $\beta: S \ra \Delta$ if and only if $S$ is isomorphic
to the compactification $\overline{\cU/L}$ of the quotient surface
$\cU/L$ for some subgroup $L$ of $\Gamma^*$ of finite index.}

{\em Proof.} Let $L$ be a subgroup of $\Gamma^*$ of finite index.
If $L \subseteq \Gamma$, there is a classical procedure that turns
the orbit space $\cU/L$ into a Riemann surface and compactifies
it. It is well-known that the canonical projection
$\overline{\cU/L} \ra \overline{\cU/\Gamma} \cong \hat{\CC}$ is a
Belyi map on the Riemann surface $S:= \overline{\cU/L}$. Now
Lemma~(2.4) implies that $S$ also admits a Belyi map when viewed
as a Klein surface. If $L \not\subseteq \Gamma$ then $L^+ : = L
\cap \Gamma$ is a subgroup of index $2$ in $L$. If we choose any
element ${\sigma}$ in $L \backslash L^+$ then ${\sigma}$ induces
an antiholomorphic involution $\tau$ on $\cU/L^+$ which can be
extended to the compactification $\overline{\cU/L^+}$ (denoted
$\tau$ again) and the following diagram commutes by Lemma~(3.1):
\[\xymatrix{\cU/L^+ \ar[r]^\tau \ar[d]^J & \cU/L^+ \ar[d]^J \\ \CC
\ar[r]^{\rm c} & \CC}.\] Hence $\overline{\cU/L} :=
\overline{\cU/L^+}/\langle\tau\rangle$ is a Klein surface and the
$J$-function induces a map $\beta: \overline{\cU/L} \ra \Delta$
which is a Belyi map because its complex double $J:
\overline{\cU/L^+} \ra {\hat{\CC}}$ is a Belyi map on the Riemann
surface
$\overline{\cU/L^+}$. \\
To prove the converse, let now $\beta: S \ra \Delta$ be a Belyi
map. If $S$ is the Klein surface associated with a Riemann surface
then $S$ also admits a Belyi map when viewed as a Riemann surface
(see Lemma~(2.4)), and, by Corollary~2 in \cite{JS2}, $S$ is
isomorphic to the compactification $\overline{\cU/L}$ for some
subgroup $L$ of finite index in $\Gamma$, in fact. Otherwise the
complex double $S^{\rm c}$ is connected and $\beta^{\rm c}: S^{\rm
c} \ra \Delta^{\rm c} = {\hat{\CC}}$ is a Belyi map. Again by
Corollary~2 in \cite{JS2}, the Riemann surface $S^{\rm c}$ is
isomorphic to the compactification $\overline{\cU/H}$ of $\cU/H$
for some subgroup $H$ of $\Gamma$ of finite index. The complex
double $S^{\rm c}$ comes with an antiholomorphic involution $\tau$
such that the following diagram commutes:
\[\xymatrix{ S^{\rm c} \ar[r]^\tau \ar[d]^{\beta^{\rm c}} & S^{\rm c} \ar[d]^{\beta^{\rm c}}
\\ {\hat{\CC}} \ar[r]^{\rm c} & {\hat{\CC}}}\]
It induces an antiholomorphic involution on $\overline{\cU/H}$
(denoted $\tau$ again) which restricts to an involution of the
open part $(\cU/H)^o$ corresponding to $(\beta^{\rm
c})^{-1}({\hat{\CC}} \backslash \{0,1,\infty\})$. Because $H$ acts
without fixed points on the preimage $\cU^o$ of $(\cU/H)^o$ the
map $\tau: (\cU/H)^o \ra (\cU/H)^o$ can be lifted to a continuous
map ${\sigma}: \cU \ra \cU$, see Part~C of Chapter~6 in
\cite{Ber}. The map ${\sigma}$ is antiholomorphic because the
canonical projection $\cU^o \ra (\cU/H)^o$ is locally
biholomorphic. Both ${\sigma}$ and the reflection $R$ induce the
complex conjugation $c$ on ${\hat{\CC}} \backslash \{0,1,\infty\}$
(by Lemma~(3.1)). Hence the composition ${\sigma} \circ R$ is a
transformation of the unramified covering $J:\cU^o \ra {\hat{\CC}}
\backslash \{0,1,\infty\}$ and therefore an element of the
covering group $\Gamma$. Thus the group $L$ generated by $H$ and
${\sigma}$ is a subgroup of finite index in the extended modular
group $\Gamma^*$ and we have $S \cong S^{\rm c}/\langle \tau
\rangle \cong \overline{\cU/H}/\langle \tau \rangle \cong
\overline{\cU/L}$ as was to be shown. \hfill $\Box$

For each positive integer $n$ let $\Gamma^*(n)$ denote the subgroup
of $\Gamma^*$ generated by the principal congruence subgroup
$\Gamma(n)$ and the reflection $R$. If is easy to see that
$\Gamma(n)$ is of index $2$ in $\Gamma^*(n)$. A fundamental region
for the action of $\Gamma^*(2)$ on the upper halfplane $\cU$ is
given by the ideal hyperbolic triangle $\cT$ with vertices $0$, $1$
and $\infty$ bounded by the imaginary axis, the line ${\rm Re}(z) =
1$ and the part of the circle $|z-\frac{1}{2}| = \frac{1}{2}$ in the
upper halfplane. The group $\Gamma^*(2)$ is generated by the
reflections $c_0$, $c_1$ and $c_\infty$ in the sides of the triangle
$\cT$, and the presentation of $\Gamma^*(2)$ is
\[\langle c_0, c_1, c_\infty \, | \, c_0^2 = c_1^2 = c_\infty^2 = 1
\rangle.\] Recall the $\lambda$-function $\lambda : \cU \ra
{\hat{\CC}} \backslash \{0,1,\infty\}$ is defined by
\[\lambda(\tau) = \frac{\wp_\tau(\frac{1+\tau}{2}) - \wp_\tau(\frac{1}{2}) }
{\wp_\tau(\frac{1+\tau}{2}) - \wp_\tau (\frac{\tau}{2})}\] where
$\wp_\tau$ denotes the Weierstra{\ss} $\wp$-function associated
with the lattice $\ZZ + \ZZ\tau$. The $\lambda$-function turns
$\cU$ into the universal covering space of ${\hat{\CC}} \backslash
\{0,1,\infty\}$ with covering group $\Gamma(2) = \langle S_0, S_1,
S_\infty \, | \, S_0 S_1 S_\infty = 1 \rangle$.

{\bf (3.3) Lemma.} {\em The $\lambda$-function induces an
isomorphism of Klein surfaces \[\cU/\Gamma^*(2)
\,\,\stackrel{\sim}{\ra} \,\, \Delta \backslash \{0,1,
\infty\}.\]}

{\em Proof.} As in Lemma~(3.1) it suffices to check that the diagram
\[\xymatrix{ \cU \ar[r]^R \ar[d]^\lambda & \cU \ar[d]^\lambda \\
{\hat{\CC}} \backslash \{0,1, \infty\} \ar[r]^{\rm c} &
{\hat{\CC}} \backslash \{0,1,\infty\}}\] commutes; this is obvious
since $\wp_{-\bar{\tau}}(\bar{z}) = \overline{\wp_\tau(z)}$ and
since $\wp_\tau$ is double periodic. \hfill $\Box$

{\bf (3.4) Theorem.} {\em A compact Klein surface $S$ admits a
Belyi map $\beta: S \ra \Delta$ if and only if $S$ is isomorphic
to the compactification $\overline{\cU/M}$ for some subgroup $M$
of $\Gamma^*(2)$ of finite index.}

{\em Proof.} We may consider any subgroup of $\Gamma^*(2)$ as a
subgroup of $\Gamma^*$. Thus one direction is an immediate
consequence of Theorem~(3.2). The other direction can be proved in
the same way as in Theorem (3.2) by replacing the $J$-function by
the $\lambda$-function and Corollary~2 in \cite{JS2} by
Corollary~1 in \cite{JS2}. (Actually both directions can be proved
this way.) \hfill $\Box$

Let $l_0$, $l_1$ and $l_\infty$ be positive integers or $\infty$. We
call the group
\[\Gamma^*(l_0, l_1, l_\infty) = \langle r_0, r_1, r_\infty \, | \,
r_0^2 = r_1^2 = r_\infty^2 = (r_0r_1)^{l_\infty} = (r_\infty
r_0)^{l_1} = (r_1 r_\infty)^{l_0} = 1 \rangle \] an {\em extended
triangular group}. It acts discontinuously and faithfully on a
Riemann surface $\cX = {\hat{\CC}}$, $\CC$ or $\cU$ as $\sum
l_i^{-1}
>1$, $=1$ or $<1$. The generators $r_i$ correspond to reflections in
the sides of a triangle $\cT$ with angles $\pi/l_i$. We view the
triangle group
\[\Gamma(l_0,l_1,l_\infty) = \langle T_0, T_1, T_\infty \, | \,
T_i^{l_i} = T_0 T_1 T_\infty = 1 \rangle\] as a subgroup of
$\Gamma^*(l_0,l_1,l_\infty)$ via $T_0, T_1, T_\infty \mapsto r_1
r_\infty, r_\infty r_0, r_0r_1$. It is easy to see that it is of
index $2$ with a representative for the nontrivial coset given by
any of the reflections $r_i$. We recall that there is a
holomorphic mapping
\[j:\cX \ra {\hat{\CC}}\]
mapping the vertices of the triangle $\cT$ to $0$, $1$ and
$\infty$ and inducing an isomorphism $\cX/\Gamma(l_0, l_1,
l_\infty) \,\, \stackrel{\sim}{\ra} \,\, {\hat{\CC}}$.

{\bf (3.5) Lemma.} {\em The map $j: \cX \ra {\hat{\CC}}$ induces
an isomorphism of Klein surfaces
\[ \cX/\Gamma^*(l_0, l_1, l_\infty) \,\, \stackrel{\sim}{\ra} \,\,
\Delta.\]}

{\rm Proof.} As in Lemma~(3.1) it suffices to check that the diagram
\[\xymatrix{\cX \ar[r]^{r_i} \ar[d]^j & \cX \ar[d]^j \\ {\hat{\CC}}
\ar[r]^{\rm c} & {\hat{\CC}}} \] commutes; this immediately
follows from the fact that $j$ maps each of the sides of the
triangle $\cT$ to the real line in $\CC \subset {\hat{\CC}}$.
\hfill $\Box$

{\bf (3.6) Theorem.} {\em A compact Klein surface $S$ admits a
Belyi map $\beta: S \ra \Delta$ if and only if $S$ is isomorphic
to the quotient surface $\cX/K$ for some subgroup $K$ of finite
index in an extended triangular group.}

{\em Proof.} Replacing the $J$-function by the $j$-function and
the references {\em Lemma~(3.1)} and {\em Corollary~2 in
\cite{JS2}} by {\em Lemma~(3.5)} and {\em Corollary~3 in
\cite{JS2}}, respectively, this can be proved in the same way as
in Theorem~(3.2). \hfill $\Box$

\bigskip

\bigskip

\section*{\S 4 Maps on Klein surfaces}

There is a well-developed theory of maps on Riemann surfaces.
Without giving precise definitions we remind the reader of the
main ideas which are described in detail in \cite{JS1}. Let $S$ be
a compact orientable surface. A topological map $\cM$ on $S$ is an
embedding of a connected graph $\cG$ into $S$ such that $S
\backslash \cG$ is a union of polygonal $2$-cells. The graphs may
have loops, multiple edges and {\em free} edges; these are edges
homeomorphic to a closed interval with only one vertex. A dart is
an ordered pair consisting of an edge together with an incident
vertex. Thus non-free edges contain two darts and the free edges
contain one dart, as illustrated in Figure 1, 2 and 3.

\begin{figure}[htbp] 
   \centering
  \includegraphics[width=2.5in]{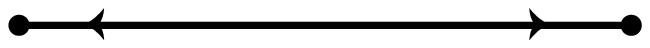}
  \caption{}
  \label{fig:example1}
\end{figure}

\begin{figure}[htbp] 
   \centering
 \includegraphics[width=2in]{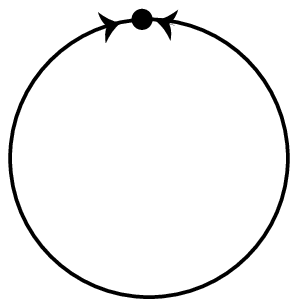}
   \caption{}
   \label{fig:example2 }
\end{figure}

\begin{figure}[htbp] 
   \centering
 \includegraphics[width=2.5in]{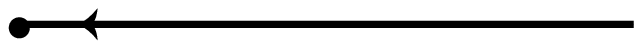}
  \caption{}
   \label{fig:example3}
\end{figure}

Let $\Omega$ denote the set of darts of $\cM$. Then there is an
involution $x$ defined on $\Omega$ as follows. If $\alpha, \beta$
are darts on a non-free edge, we define $x(\alpha) = \beta$. If
$\alpha$ is a dart on a free edge we define $x(\alpha) = \alpha$.
Thus we have a natural 1-1 correspondence between the set of
orbits of $\langle x \rangle$ on $\Omega$ and the set of edges of
$\cM$. Define a permutation $y$ of $\Omega$ that cyclically
permutes the darts of $\Omega$ that are incident with a vertex $v$
of $\cG$. Precisely, $y$ takes a dart incident with $v$ and
transfers it, using the orientation of $S$ to the next dart
incident with $v$ in an anticlockwise direction. Thus we have a
natural 1-1 correspondence between the set of orbits of $\langle y
\rangle$ on $\Omega$ and the set of vertices of $\cM$. It is then
easy to see that $y^{-1}x$ cyclically permutes the darts around a
face in an anticlockwise direction, so that we have the relations
$x^2 = y^m = (y^{-1}x)^n = 1$, where $m$ is the least common
multiple of the vertex valencies and $n$ is the least common
multiple of the face sizes. Let $G$ denote the subgroup of the
permutation group of $\Omega$ generated by $x$ and $y$. It acts
transitively on $\Omega$ as $\cG$ is connected.

This motivates us to consider algebraic maps. An {\em algebraic
map} is a quadruple $(G, \Omega, x, y)$ where $G$ is a group
acting transitively and faithfully on a finite set $\Omega$ and
$x$, $y$ are generators of $G$ with $x^2 = 1$. Thus associated
with a topological map $\cM$ we have an algebraic map ${\rm Alg}
\cM$. In \cite{JS1}, we showed how to invert this association and
pass from an algebraic to a topological map. We do this by
considering the universal topological map $\hat{\cM}(m,n)$ of type
$(m,n)$. This is the regular map on one of the three simply
connected surfaces $\cX$ all of whose vertices have valency $m$
and all of whose faces have size $n$. For example, if $m=n=4$, we
have the chessboard tessellation of the Euclidean plane. The
associated algebraic map ${\rm Alg}\hat{\cM}(m,n)$ is $(\Gamma,
\hat{\Omega}, X, Y)$ where $\Gamma$ is the triangle group
\[\Gamma(2,m,n) = \langle X, Y \, | \, X^2 = Y^m = (Y^{-1}X)^n = 1
\rangle\] acting on $\hat{\Omega} = \Gamma$ via right
multiplication. Given an algebraic map $(G, \Omega, x, y)$ of type
$(m,n)$ there is an obvious epimorphism $\Gamma \ra G$ and thus
$\Gamma$ acts transitively on $\Omega$. Let $M$ be the stabilizer
of some element in $\Omega$. Then $M$ is a Fuchsian group and so
$\cX/M$ is a Riemann surface. Thus with the algebraic map
$(G,\Omega, x, y)$ we have associated the Riemann surface $\cX/M$.
If $(G,\Omega, x, y)$ is the algebraic map associated with a
topological map $\cM$ on the surface $S$, then $\cX/M$ is
homeomorphic to $S$; in particular we have turned $S$ into a
Riemann surface. We recall that $\cX/\Gamma \cong {\hat{\CC}}$ and
that the natural projection $\cX/M \ra \cX/\Gamma \cong
{\hat{\CC}}$ is a Belyi function on the Riemann surface $\cX/M$.

Vice versa, we can use Belyi functions to construct maps. On
${\hat{\CC}}$ we define the trivial map with one dart with a
vertex $0$ and free edge along the equator with second endpoint at
$1$ and a face center at $\infty$. If $\beta: X \ra {\hat{\CC}}$
is a Belyi function on the compact Riemann surface $X$ then we
define a map $\cM$ on $X$ by using $\beta$ to pull back the
trivial map on ${\hat{\CC}}$. Then the Riemann surface $\cX/M$
constructed from $\cM$ as above is isomorphic to $X$ (see
Section~6 in \cite{JS2}).

A theory of maps on (possibly non-orientable) compact surfaces
with boundary was developed in \cite{BS}. Also see \cite{JS3}.
This was based on the work of \cite{JS1} described above. As
above, a topological map $\cM$ on a compact surface $S$ with
boundary is an embedding of a connected graph $\cG$ as above into
$S$ with fairly obvious conditions. For instance, an edge of $\cG$
not lying on the boundary $\partial S$ can intersect $\partial S$
in at most two points and these cannot be interior points of the
edge. At a vertex $p$ on $\partial S$ there may 0, 1 or 2 edges
emanating from $p$ on $\partial S$. Each face, i.e.\ connected
component of $S \backslash \cG$, is homeomorphic to the disc $D=
\{ z \in \CC \, | \, |z| < 1\}$ or to the half disc $E = \{z \in
\CC \, | \, |z| < 1, {\rm Im}(z) \ge 0\}$. First we define a blade
to be a half-dart. So an interior dart will have two blades,
$b_1$, $b_2$, one which points above the edge, and one which
points below the edge. If the dart lies on $\partial S$, then the
dart will have only one blade. (For a rigorous account see
\cite{BS}.) We can define three involutions $\tau, \lambda, \rho$
on the set of blades. At a vertex in the interior $S^o$, the
involution $\tau$ transposes the upper and lower blade  and it
fixes any blade at a vertex on $\partial S$. For instance
$\tau(b_1) = b_2$ and $\tau(b_3) = b_4$ in Figure~4.

\begin{figure}[htbp] 
   \centering
  \includegraphics[width=2in]{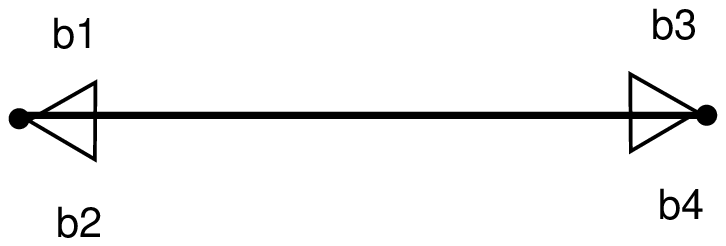}
   \caption{}
   \label{fig:example4}
\end{figure}

We call $\tau$ the {\em transverse reflection}. The so-called {\em
longitudinal reflection} $\lambda$ usually interchanges the two
upper `upper' blades and also the two `lower' blades of an edge.
For instance $\lambda(b_1) = b_3$ and $\lambda(b_2) = b_4$ in
Figure~4. If the edge intersects the boundary, or if it is a free
edge, then $\lambda$ may fix a blade. The last involution we
define is the {\em rotary reflection} $\rho$. As a blade $b_1$
points in a particular direction (either clockwise or
anti-clockwise), we follow this direction and if we meet another
blade, $b_5$, we define $\rho(b_1) = b_5$. As blade $b_5$ points
in the opposite direction, $\rho(b_5) = b_1$ as in Figure 6.

\begin{figure}[htbp] 
   \centering
   \includegraphics[width=2in]{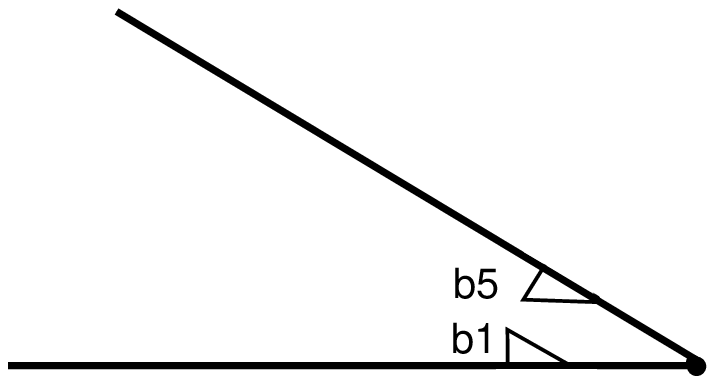}
   \caption{}
   \label{fig:example6}
\end{figure}

Again, $\rho$ may fix a blade which corresponds to boundary
behaviour.

Notice that $\tau$ and $\rho$ generate a dihedral group that fixes
a vertex, $\rho$ and $\lambda$ generate a dihedral group that
fixes a face and that $\lambda$ and $\tau$ generate a dihedral
group fixing an edge. The group $G^*$ generated by $\tau$,
$\lambda$ and $\rho$ obeys the relations
\[\tau^2 = \lambda^2 = \rho^2 = (\tau\lambda)^2 = (\tau \rho)^m =
(\rho \lambda)^n =1\] where $m$ is the least common multiple of
the vertex valencies and $n$ is the least common multiple of the
face valencies. In particular there is an obvious epimorphism
$\phi: \Gamma^*(2,m,n) \ra G^*$. Let $K \subseteq \Gamma^*(2,m,n)$
denote the stabilizer of a blade. We call $K$ the {\em map
subgroup for} $\cM$. Then $K$ is of finite index in
$\Gamma^*(2,m,n)$ and the orbit space $\cX/K$ is homeomorphic to
$S$. By Theorem~(3.6), $\cX/K$ is a Klein surface and the
canonical projection $S= \cX/K \ra \cX/\Gamma^*(2,m,n) \cong
\Delta$ is a Belyi function.

To invert this process, we define the {\em trivial map} on
$\Delta$: it consists of a single blade; regard $\Delta$ as being
a hemisphere with boundary the great circle through $0, 1,
\infty$; then the single blade consists of half a free edge with
endpoints a vertex at 0 and a free end point at 1. If $S$ is a
Klein surface and $\beta: S \ra \Delta$ is a Belyi map on $S$,
then we can pull back the trivial map on $\Delta$ to a map $\cM$
on $S$. Let $K$ denote the map subgroup for $\cM$. Then we can
find an isomorphism between the Klein surfaces $\cX/K$ and $S$
such that the following diagram commutes:
\[\xymatrix{ \cX/K \ar[r]^\sim \ar[d]^{\rm can} & S
\ar[d]^\beta \\ \cX/\Gamma^*(2,m,n) \ar@{=}[r] & \Delta }\]
Summarizing we have proved:

{\bf (4.1) Theorem.} {\em Let $S$ be a compact connected Klein
surface. Then the following statements are equivalent:

(a) $S$ admits a Belyi map $\beta: S \ra \Delta$.

(b) $S$ carries a map $\cM$ such that $S$ is isomorphic to the
quotient surface $\cX/K$ where $K$ is the map subgroup for $\cM$.

(c) $S$ is isomorphic to the quotient surface $\cX/K$ for some
subgroup $K$ of finite index in an extended triangle group of the
form $\Gamma^*(2,m,n)$.} \hfill $\Box$

School of Mathematics, University of Southampton, Highfield,
Southampton SO17 1BJ, United Kingdom.\\
{\em E-mail:} B.Koeck@soton.ac.uk, ds@maths.soton.ac.uk

\end{document}